# Technological schemes and control methods in the reconstruction of parallel gas pipeline systems under non-stationary conditions

Ilgar Aliyev

Head of Department, Azerbaijan Architecture and Construction University, Baku, Azerbaijan

## *ABSTRACT*

*The study explores technological schemes and control methods for reconstructing parallel gas pipeline systems under non-stationary conditions. It focuses on improving safety, reliability, and efficiency by automating valve control and integrating IoT-based monitoring. Automated Shut-Off Valves—Detect sudden pressure drops and block leaks instantly. Uses pressure drop rates and valve position sensors to detect and isolate leaks. Wireless pressure sensors and real-time monitoring allow remote control of gas flow. Mathematical modeling analyzes pressure variations along damaged sections to determine valve placement and activation timing. Machine learning algorithms are used in the control center to predict and verify leak locations based on sensor data. For ensuring uninterrupted gas supply in the event of an accident, this study has developed an empirical formula that determines the location of connecting pipes—activated between parallel pipelines—based on real-time pressure sensor data. The purpose of the study is to reduce gas leaks and environmental hazards, ensure an uninterrupted gas supply during emergencies, enhance decision-making through automated systems, minimize gas losses, and reduce maintenance costs.*

## 1. INTRODUCTION

The growing complexity and aging of parallel gas pipeline networks present major obstacles for contemporary energy infrastructure. As time progresses, factors like physical deterioration, outdated components, and shifting technological needs have reduced the operational dependability and safety of these systems. Consequently, there is an urgent need for reconstruction and modernization that not only updates obsolete components but also incorporates sophisticated control systems to address the variable nature of gas flow and reduce risks like gas leaks. The motivation for this research stems from the dual imperative of ensuring uninterrupted gas supply to consumers while protecting the environment from the adverse effects of uncontrolled gas leaks. Traditional systems that rely solely on pressure drop rates for valve activation have shown limitations, particularly under emergency conditions. Therefore, our approach enhances reliability by combining pressure monitoring with real-time valve position sensors and integrating Internet of Things (IoT) technologies. This dual-criteria system enables the automatic, simultaneous closure of valves on both sides of a leak, thereby isolating the compromised section of the pipeline more effectively.

The paper is structured to provide a comprehensive exploration of the technological schemes and control methods necessary for the reconstruction of parallel gas pipeline systems under non-stationary conditions [16]:

Background and Motivation: The introduction outlines the current challenges in maintaining and modernizing aging gas pipeline networks. It discusses the limitations of conventional control systems and the environmental and operational risks associated with gas leaks.

Technological Schemes and Control Methods: This section describes the proposed methodological framework, detailing the placement of automatic shut-off valves, the role of valve position sensors, and the integration of IoT and wireless communication systems. It also explains how these innovations improve the response to emergencies and maintain a stable gas supply.

Mathematical Modeling and Experimental Analysis: A detailed analysis of non-stationary gas flow in damaged pipeline sections is presented. The section introduces mathematical models, experimental data, and pressure variation profiles that help in understanding the dynamic behavior of gas flows during valve closure and subsequent system recovery.Algorithm Development and System Integration: Building on the experimental insights, this part explains the derivation of an empirical algorithm based on machine learning. The algorithm uses sensor data to determine the optimal activation times and locations for the connecting valves, ensuring minimal disruption to gas supply and enhanced safety. Conclusions and Recommendations: The final section summarizes the findings, emphasizing the practical benefits of the proposed technological schemes for modernizing existing gas pipelines. It also outlines potential directions for future research and implementation strategies. The primary goal of this research is to develop and validate an advanced control system for parallel gas pipeline networks that improves safety, operational efficiency, and reliability under non-stationary flow conditions. Specifically, the study aims to enhance emergency response: by integrating valve position monitoring sensors with traditional pressure drop criteria, the system is designed to promptly and simultaneously close valves on both sides of a leak, thereby reducing gas loss and environmental impact. Optimize Gas Flow Dynamics: Through mathematical modeling and experimental analysis, the research examines the behavior of non-stationary gas flows in damaged pipeline sections, providing insights into the pressure variations and identifying optimal control strategies.

Facilitate Modernization Using IoT: The proposed scheme incorporates IoT technology for real-time data collection, wireless communication, and cloud-based analysis, enabling centralized control and improved decision-making processes.

Ensure Continuous Gas Supply: By accurately determining the activation times and locations for pneumatic valves on connecting pipes, the system aims to maintain an uninterrupted gas supply to consumers, even in emergency scenarios [13].

The problems addressed in this study include the identification of leakage points in linear pipeline networks, the challenge of managing non-stationary gas flows during emergency valve operations, and the need to integrate modern sensor technologies into existing gas distribution systems [11,22,26].

In conclusion, this study presents a strong framework for the reconstruction of parallel gas pipeline systems by integrating advanced control techniques, sensor technologies, and IoT. The proposed solution not only enhances operational reliability and safety but also provides a scalable, cost-efficient approach to modernizing aging gas infrastructure.

## 2. MATERIALS AND METHODOLOGY

The methodological foundation of this study lies in the development of technological frameworks and algorithms aimed at enhancing automated control systems for non-stationary gas transportation processes during the reconstruction of parallel gas pipelines. This can address key issues such as increasing the efficiency, reliability, and safety of gas transportation systems at the required level.

Thus, for the reconstruction of existing parallel gas pipelines, the use of modern separation devices and elements is proposed, along with the reliable management of non-stationary gas flows. The following structural scheme is proposed for this purpose (Figure 1).

In the case of an emergency situation in parallel pipelines, the diameter of the connecting pipe, which is installed to transfer the gas flow from one pipeline configuration to another, is taken to be equal to the diameter of the gas pipeline. At the same time, from the moment the gas pipeline is put into operation, along with its initial parameters ($P_1$, $P_2$, $G_0$, L, d), the distance between the connecting pipes ($\ell$) is also recorded in the technical passport of the gas pipeline.

According to Figure 1, the installation of automatic shut-off valves on the parallel gas pipeline is planned as follows.

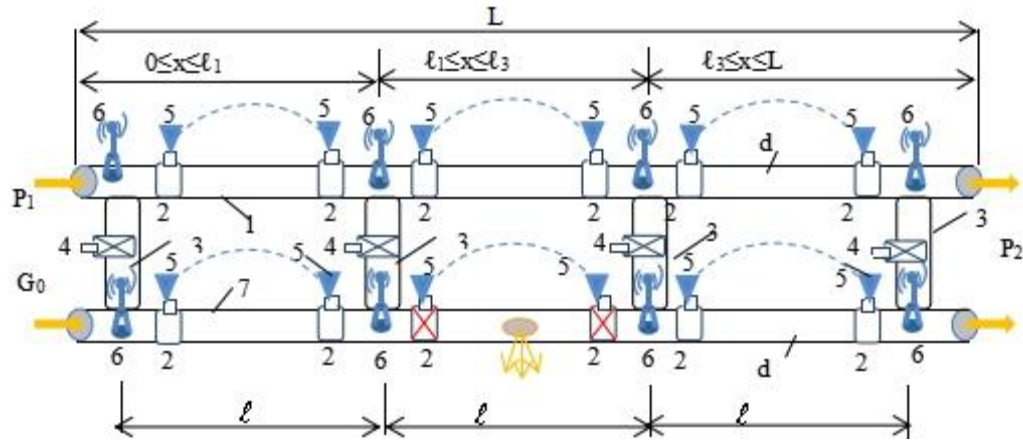

Figure 1. Proposed structural diagram for the reconstruction phase of parallel gas pipelines.
1- undamaged pipeline; 2- pipeline isolation valves; 3- pipeline connecting the parallel pipelines (connecting pipes); 4- isolation valves of the connecting pipes; 5- valve position monitoring sensor; 6- wireless pressure sensor; 7- damaged pipeline.

L - length of the gas pipeline, m; $P_1$ - pressure at the beginning of the gas pipeline in stationary mode, Pa; $P_2$ - pressure at the end of the gas pipeline in stationary mode, Pa; $G_0$ - mass flow rate at the beginning of the gas pipeline, $\frac{Pa \cdot \sec}{m}$; d - diameter of the pipeline, m; $\ell$ - distance between the connecting pipes, m.

They are to be installed after the connecting pipe at the beginning of the gas pipeline and before the connecting pipe at the end of the gas pipeline. In the other sections of the parallel gas pipeline, automatic shut-off valves are to be installed on both the right and left sides of each pipeline in the network.

For the efficient management of gas pipelines in emergency situations, the automatic activation of valves has proven to be effective in practice. Previously, valves calibrated based on the rate of pressure drop in the gas pipeline would automatically activate each valve separately if the pressure drop rate exceeded 0.1 MPa/min during emergency conditions [8]. Recently, however, with the application of valve position monitoring sensors, the automatic activation of valves has been successfully tested in practice and has been confirmed as a reliable technology [12].

The main advantage of the scheme we propose is that, in the event of a gas leak, the calibration (adjustment) of the valves is based not only on the rate of pressure drop but also on the activation of the valve position monitoring sensor. This improves the reliability of valve operation in emergency conditions. Specifically, when the integrity of one of the parallel gas pipeline pipes is compromised, the valve of the nearest shut-off valve to the leak location will close when the pressure drop reaches the calibrated value (the pressure drop rate reaches 0.2-0.5 MPa per minute [5]). The sensor monitoring the position of the closed valve automatically activates, which also activates the sensor monitoring the valve on the other side of the leak. The activated sensor automatically opens the valve it is monitoring and closes the shut-off valve. Therefore, both the rate of pressure drop and the activation of valve sensors will simultaneously trigger and close the valves on both sides of the leak. As shown in Figure 1, each of the automatic shut-off valves installed on the parallel gas pipeline is equipped with a position monitoring sensor. The sensors of the automatic shut-off valves installed on the pipeline configurations between the connecting pipes communicate with each other through the technical operation capabilities of the Internet of Things gateway. This allows the activation of one sensor to automatically trigger the sensor it is connected to. If the automatic shut-off valve on the left side of the leak point closes due to the limited pressure drop rate (Figure 1), the sensors designed to activate the valves will automatically close the valve on the right side of the leak. In other words, when the valve on the left side of the leak is automatically closed due to the pressure drop rate, the sensor monitoring the position of that valve activates, which triggers the sensor on the right side of the leak, connected via wireless communication channels. This will automatically close the valve on the right side of the leak. Thus, the automatic shut-off valves on both the right and left sides of the leak point will be closed simultaneously.

For the remote activation of pneumatic valves installed on the connecting pipes of the parallel gas pipeline, wireless pressure sensors are proposed [3,14]. As shown in Figure 1, it is suggested that the pressure sensors be installed on the section of the pipeline where the connecting pipes of the gas pipelines are located. These sensors allow for the monitoring of pressure variations in the pipeline configurations, enabling the control center to track the data via the internet.

It is known that the following technologies are available for predicting the technical operation of pipelines and managing emergency situations: the Internet of Things (IoT) [6]; machine learning and artificial intelligence [1,20]; advanced communication protocols such as 5G for real-time data transmission [18,24]; Arduino Uno controller - cloud storage for data collection and analysis [7,15,27]; hybrid models [4]; real-time monitoring or artificial intelligence [2,18]. These technologies are playing an increasingly important role in complex gas transportation systems and should be considered for the efficient operation of gas pipelines. It is well known that trunk and distribution gas pipelines are mainly designed in remote and inaccessible areas of residential settlements. Therefore, previously, due to the difficulty of monitoring the technological regime of these gas pipelines, their operational costs were correspondingly high. However, now, with the application of innovation, changes occurring in the gas pipelines can be monitored and tracked remotely via the internet through devices installed on them. As an example, the application of the Internet of Things (IoT) technology for monitoring and managing the non-stationary operating modes of the linear section of the gas pipeline can be highlighted [26]. The Internet of Things (IoT) technology is a system that enables devices containing hardware, software, and networks to interact, cooperate, and share information. Currently, IoT technology sensors remotely monitor the technological and emergency modes of these pipelines in real-time via the internet and send the collected data directly to control centers. Thus, if a leak is detected, the system automatically activates (closes) the valves to isolate the damaged section of the pipeline from the main section. Numerous local and international studies have addressed these issues [9,10,19,21,23,25,29].

In these works, the authors describe how the integration of the mentioned technologies and methodologies into gas pipelines will enhance the efficiency, safety, and reliability of gas pipeline systems. The materials presented, however, do not extensively address the methods for detecting leakage points in the linear pipeline or the analysis of non-stationary gas flow parameters along the entire length of the network as a result of this process. The key difference in our research is that it analyzes various transition processes of non-stationary gas flow in parallel gas pipelines. Specifically, in the event of a breach in the integrity of one of the pipeline sections of parallel gas pipelines, the process of activating the valves located to the left and right of the damaged section, as well as the remote activation of pneumatic valves installed on the connecting pipes for uninterrupted gas supply to consumers, are analyzed together, as shown in Figure 1. From the moment the valves on the pipeline are closed ($t=t_1$), the analysis includes the non-stationary gas flow parameters in three different sections of the damaged pipeline until the activation of the connecting pipes (change in pressure values) ($t=t_2$). Additionally, the technological scheme for accepting the control method of gas dynamic processes of these transition processes through a control center is developed. In this research, the use of pressure sensors for monitoring the transition processes, alongside the valve position sensor, is proposed for the localization of the damaged section in the emergency condition of parallel gas pipelines. The methodology of the proposed control system for the parallel gas pipeline is based on the following scheme (Figure 2).

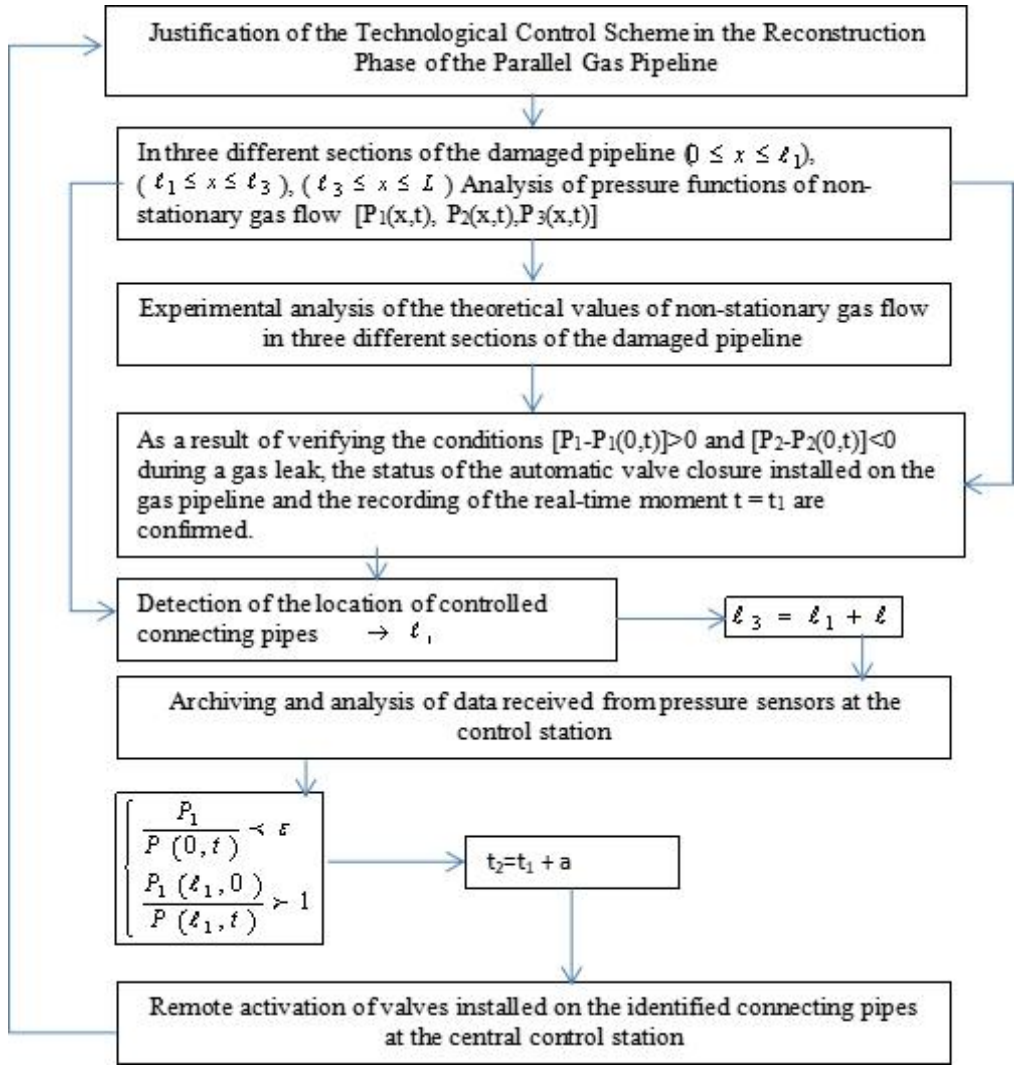

Figure 2. Methodological diagram of the parallel gas pipeline control system.

The analysis of the methodological scheme of the control system for the parallel gas pipeline (Figure 2) shows that data collection through the Internet of Things (IoT), as well as data storage and analysis using cloud storage networks, is proposed to be integrated into a control center, and a reporting scheme is suggested for an efficient wireless control system method. Therefore, when the integrity of any pipeline in the parallel gas pipeline is compromised, the valves located on either side of the leak are closed at time t=t1. This occurs due to the pressure drop rate reaching its threshold value and the activation of the sensors installed on the valves, effectively isolating the damaged section from the main structure of the parallel gas pipeline system. Based on Figure 1, the mathematical models of the gas dynamic process occurring in three different sections of the damaged pipeline ( $0 \leq x \leq \ell_1$, $\ell_1 \leq x \leq \ell_3$,

$\ell_3 \leq x \leq L$ ) from the moment the valves are closed have been thoroughly analyzed in the study [17] we investigated, and the pressure functions of the non-stationary gas flow for the three sections are represented by the following expressions [17].

$$0 \leq x \leq \ell_1$$

$$P_1(x,t) = P_1(x,t_1) + \frac{c^2}{\ell_1}\int_{t_1}^{t} G_0(\tau)d\tau + \frac{2c^2}{\ell_1}\sum_{n=1}^{\infty}\cos\frac{\pi n x}{\ell_1}\int_{t_1}^{t} G_0(\tau)\cdot e^{-\alpha_3(t-\tau)}d\tau +$$
$$+\frac{c^2}{2a\ell_1}\frac{dP_1(0,t_1)}{dx}\int_{t_1}^{t}\left[1+2\sum_{n=1}^{\infty}\cos\frac{\pi n x}{\ell_1}e^{-\alpha_3\tau}\right]d\tau \quad (1)$$

$$\ell_1 \leq x \leq \ell_3$$

$$P_2(x,t) = P_2(x,t_1) - \frac{c^2}{(\ell_3-\ell_1)}\int_{t_1}^{t} G_{ut}(\tau)\left[1+2\sum_{n=1}^{\infty}(-1)^n\cos\frac{\pi n(\ell_2-\ell_3)}{\ell_3-\ell_1}\right]\times$$
$$\times\cos\frac{\pi n(x-\ell_1)}{\ell_3-\ell_1}e^{-\alpha_4(t-\tau)}d\tau \quad (2)$$

$$\ell_3 \leq x \leq L$$

$$P_3(x,t) = P_3(x,t_1) - \frac{c^2}{(L-\ell_3)}\int_{t_1}^{t} G_s(\tau)\left[1+2\sum_{n=1}^{\infty}(-1)^n\cos\frac{\pi n(x-\ell_3)}{L-\ell_3}\cdot e^{-\alpha_5(t-\tau)}\right]d\tau -$$
$$-\frac{c^2}{2a(L-\ell_3)}\frac{dP_3(L,t_1)}{dx}\int_{t_1}^{t}\left[1+2\sum_{n=1}^{\infty}(-1)^n\cos\frac{\pi n(x-\ell_3)}{L-\ell_3}\cdot e^{-\alpha_5\tau}\right]d\tau \quad (3)$$

Here, $\alpha_3 = \frac{\pi^2 n^2 c^2}{2a\ell_1^2}$, $\alpha_4 = \frac{\pi^2 n^2 c^2}{2a(\ell_3-\ell_1)^2}$, $\alpha_5 = \frac{\pi^2 n^2 c^2}{2a(L-\ell_3)^2}$

## 3. EXPERIMENTAL ANALYSIS OF THE THEORETICAL VALUES OF NON-STATIONARY GAS FLOW

Let us assume that for the purpose of localizing the gas leak, the valves installed on the pipelines are automatically closed via the valves at $t=t_1=300$ seconds. Then, for each of the three different sections of the parallel gas pipeline, we accept the following given parameters of the specific gas pipeline in order to determine the regularity of the change in non-stationary gas flow for the considered case [17].

$P_1(0,t_1) = 13.36 \cdot 10^4 \text{ Pa}$; $P_1(5000,t_1) = 12.82 \cdot 10^4 \text{ Pa}$;
$P_1(10000,t_1) = P_2(10000,t_1) = 12.19 \cdot 10^4 \text{ Pa}$; $P_2(14500,t_1) = 11.56 \cdot 10^4 \text{ Pa}$;

$P_2(20000, t_1) = P_3(20000, t_1) = 11.24 \cdot 10^4 \, Pa; \quad G_0 = 10 \, \dfrac{Pa \cdot san}{m}$

$P_3(25000, t_1) = 10.86 \cdot 10^4 \, Pa; \; P_2 = 11 \cdot 10^4 \, Pa; \; P_3(30000, t_1) = 10.4 \cdot 10^4 \, Pa; \; 2a = 0.1 \dfrac{1}{sec}$

$P_1 = 14 \cdot 10^4 \, Pa; \; L = 3 \cdot 10^4 \, m; \; \ell_1 = 1 \cdot 10^4 \, m; \; \ell_2 = 1.45 \cdot 10^4 \, m; \; \ell_3 = 2 \cdot 10^4 \, m; \; c = 383.3 \dfrac{m}{sec};$

For the considered period, it can be assumed that $G_0(t)$=constant [17]. Then, using equation (1), we calculate the pressure values of the non-stationary gas flow along the length of the damaged pipeline section 1 every 120 seconds, and record them in the table below (Table 1).

Table 1. Pressure variation values along the length ($0 \leq x \leq \ell_1$) of the damaged gas pipeline section 1 as a function of time.

| x. km | $P_1(x.t)$, 10⁴ Pa | | | | | |
|---|---|---|---|---|---|---|
| | t=t₁+a. sec | | | | | |
| | a=0 | a=120 | a=240 | a=360 | a=480 | a=600 |
| 0 | 13.36 | 14.58 | 15.46 | 16.35 | 17.24 | 18.13 |
| 5 | 12.82 | 13.67 | 14.55 | 15.44 | 16.33 | 17.22 |
| 10 | 12.19 | 12.91 | 13.8 | 14.69 | 15.57 | 16.46 |

It can be seen from Table 1 that after the automatic valves on the damaged pipe are closed, the pressure along the length of section 1 of the gas pipeline increases over time. Specifically, at x = 0, during the period from a = 0 to 600 seconds, the pressure increases from 13.36 × 10⁴ Pa to 18.13 × 10⁴ Pa. At x = 5 km, the pressure at a = 0 seconds is 12.82 × 10⁴ Pa, and after 600 seconds, it increases to 17.22 × 10⁴ Pa. Similarly, at the final point x = 10 km, the pressure at a = 0 seconds is 12.19 × 10⁴ Pa, and at a = 600 seconds, it increases to 16.46 × 10⁴ Pa.
That is, during the period from a = 0 to 600 seconds, the rate of increase along the length is 1.35. Accordingly, the pressure values of the non-stationary gas flow along the length of the second section of the damaged pipeline are calculated every 120 seconds and recorded in the following table (Table 2).

Table 2. Pressure variation along the length and with respect to time in the second section ($\ell_1 \leq x \leq \ell_3$) of the damaged gas pipeline.

| x. km | $P_2(x.t)$, 10⁴ Pa | | | | | |
|---|---|---|---|---|---|---|
| | t=t₁+a, sec | | | | | |
| | a=0 | a=120 | a=240 | a=360 | a=480 | a=600 |
| 10 | 12.19 | 10.59 | 9.7 | 8.81 | 7.93 | 7.04 |
| 14.5 | 11.56 | 10.52 | 9.63 | 8.74 | 7.86 | 6.97 |
| 20 | 11.24 | 10.01 | 9.13 | 8.24 | 7.35 | 6.46 |

From Table 2, it is observed that in the second section of the damaged pipeline, the pressure decreases as you move from the starting point to the endpoint. Specifically, at the starting point (x=10 km), the pressure decrease rate is 1.73 over a period of 600 seconds. At x=14.5 km, the rate of decrease is 1.65, and at the final point (x=20 km), the decrease rate is again 1.73. Based on the analysis, we can conclude that the pressure decrease rates at the starting and ending points of the second section are equal, while at the leak point, the rate of decrease increases.
Following this, the pressure values along the length of the third section of the damaged pipeline will be calculated and noted in Table 3 every 120 seconds.

Table3. Pressure variation values along the length ($\ell_1 \leq x \leq \ell_3$) of the third section of the damaged gas pipeline and over time.

| x. km | $P_3(x,t)$, $10^4$ Pa | | | | | |
|---|---|---|---|---|---|---|
| | $t=t_1+a$. sec | | | | | |
| | a=0 | a=120 | a=240 | a=360 | a=480 | a=600 |
| 20 | 11.24 | 10.52 | 9.63 | 8.74 | 7.86 | 6.97 |
| 25 | 10.86 | 10.01 | 9.13 | 8.24 | 7.35 | 6.46 |
| 30 | 10.4 | 9.19 | 8.3 | 7.41 | 6.52 | 5.63 |

From Table 3, it can be seen that in the third section, the pressure decreases along the direction from the starting point to the end point over time. For instance, at the starting point (x = 20 km), at t = 0 seconds, the pressure is 11.24×10⁴ Pa. At t = 240 seconds, it decreases to 9.52×10⁴ Pa, and at t = 600 seconds, it further decreases to 6.97×10⁴ Pa. Thus, the rate of decrease at t = 600 seconds is 1.6. At the x = 25 km point, at t = 0 seconds, the pressure is 10.86×10⁴ Pa, and at t = 600 seconds, the pressure decreases to 6.46×10⁴ Pa, giving a rate of decrease of 1.68. At the x = 30 km endpoint, at t = 600 seconds, the pressure decreases to 5.63×10⁴ Pa, resulting in a decrease rate of 1.85. Based on the analysis, we can conclude that the rate of pressure decrease increases from the starting point to the end point of the third section. As a result of the sharp decrease in pressure in the final section of the damaged gas pipeline, the volume of gas supplied to consumers will also decrease.

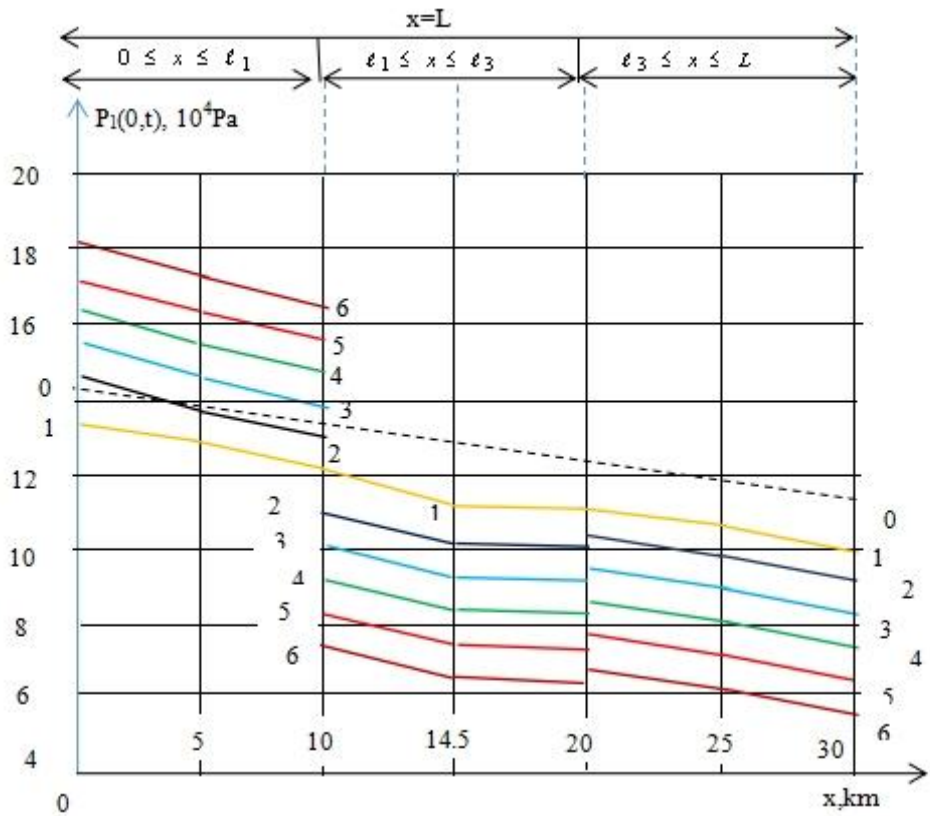

Figure 3. The graph of pressure variation over time along the length of the gas pipeline after the valves at
x=ℓ1 and x=ℓ3 are closed at t=t1. (0 - t=0 stationary state, 1 - t=t1, 2 - t=t1+120 sec, 3 - t=t1+240 sec, 4 - t=t1+360 sec, 5 - t=t1+480 sec, 6 - t=t1+600 sec).

Using Tables 1, 2, and 3, we can construct a graph (Figure 3) illustrating the distribution of pressure over time along the length of the damaged pipeline in the 1st, 2nd, and 3rd sections of the complex gas pipeline.

As shown in Figure 3, in the case where the integrity of the parallel gas pipeline arrangements is not violated, i.e., in a stationary regime, the pressure values along the length of both pipelines are equal and do not change over time (0-0 line). When the tightness of one of the pipelines is compromised (according to Figure 3, the 7th arrangement), no change occurs in the gas flow dynamics in the first pipeline, meaning it continues to operate in a stationary regime (0-0 line). In the pipeline with compromised integrity, however, the pressure values along the length decrease over time (as shown by the 1-1 line in Figure 3). At $t=t_1=300$ seconds, the pressure at the starting point x=0 changes from $P_1=14\times10^4$ Pa to $13.36\times10^4$ Pa, at x=10 km it decreases from $P(\ell_1,0)=13\times10^4$ Pa to $12.19\times10^4$ Pa, at the leak point (x=15 km) from $P(\ell_2,0)=12.6\times10^4$ Pa to $11.6\times10^4$ Pa, at x=20 km from $P(\ell_3,0)=12\times10^4$ Pa to $11.24\times10^4$ Pa, and at the final point x=L the pressure decreases from $P_2=11\times10^4$ Pa to $10.4\times10^4$ Pa.

Based on the analysis, it can be concluded that the pressure decreases along the length of the damaged pipeline, and the rate of decrease at the leak point is higher compared to the other sections.

As noted, in the case of an accident, at t=t1, the valves in the sections x=$\ell_1$ and x=$\ell_3$ of the pipeline are automatically closed.

## 4. CREATİON OF TECHNOLOGİCAL MANAGEMENT SYSTEMS AND ALGORİTHMS.

Based on the theoretical analysis of the gas pipeline's non-stationary flow regimes, it can be concluded that the opening of pneumatic valves on the connecting pipes at t=t1 in the damaged sections is not advisable according to the gas dynamic process requirements. This is because, in the damaged section at x=$\ell_1$, the pressure is lower ($12.19\times10^4$ Pa $<13\times10^4$ Pa) compared to the pressure in the undamaged pipeline. As a result, the gas flow will occur in the reverse direction: gas will flow from the undamaged section into the damaged section, which will disrupt the continuity and stability of the gas supply to consumers.

Therefore, it is more appropriate to wait until the connecting pipes are opened at a later time t=t2, ensuring that the pressure in the damaged section exceeds the pressure in the undamaged section. According to Figure 3, during this time, the increase in pressure at the beginning of the damaged section at x=0 surpasses the pressure increase at x=$\ell_1$. This period is crucial for ensuring the protection of the main equipment of the gas pipeline, particularly the compressor stations. To avoid potential damage, the compression ratio ε\varepsilonε of the compressor's driving aggregates must be limited, with the value of ε≤1.3 being the accepted standard to prevent over compression, represented as [$P_1/P(0,t_2)\leq1.3$] [17].

On the other hand, as shown in Figure 3, after the valves installed on the pipeline are automatically closed, the volume of gas supplied to consumers gradually decreases over time. To ensure continuous supply of the required amount of gas to consumers, it is essential to remotely activate the valves installed on the connecting pipeline sections at x=$\ell_1$ and x=$\ell_3$. For this purpose, the control center must first receive information about the activation of the valves and quickly determine their location within the pipeline.

The experimental analysis of the theoretical values of the non-stationary gas flow pressure shows that, after the valves are closed, the pressure in the starting section of the pipeline increases over time. However, in the final section, the pressure decreases with time. The analysis indicates that this particular gas dynamic behavior can only occur in the case of valve closure. Therefore, by analyzing data from the pressure sensors installed at the start and end points of the gas pipeline at the control center, this characteristic can be identified, which will signal the valve closure process.

To determine the location of these valves, we use the equation from (1) to express the pressure change at the starting point of the pipeline as follows:
The experimental analysis of the theoretical values of the non-stationary gas flow pressure indicates that, after the valves are closed, the pressure in the initial section of the pipeline increases over time. However, in the final section, the pressure decreases with time. This behavior suggests that this particular gas dynamic phenomenon occurs only in the case of valve closure. Therefore, by analyzing the data from the pressure sensors installed at the start and end points of the pipeline at the control center, this characteristic can be detected, which will indicate the valve closure process.
To determine the location of these valves, we use the equation from (1) to express the pressure variation at the starting point of the pipeline as follows:

$$P_1(0,t) = P_1(0,t_1) - \frac{2a\ell_1}{3}G_0 + \frac{2c^2}{\ell_1}G_0 \sum_{n=1}^{\infty} \frac{e^{-\alpha n^2 (t-t_1)} + e^{-\alpha n^2 t} - e^{-\alpha n^2 t_1}}{\alpha n^2} \qquad (4)$$

Here, $\alpha = \frac{\pi^2 c^2}{2a\ell_1^2}$

The simplified form of equation (4) will be as follows:

$$P_1(0,t) = P_1(0,t_1) + \frac{2c^2 G_0}{\ell_1}(t-t_1)(1-C) - 2aG_0\ell_1 \left(\frac{1}{3} + \frac{2}{\pi^2}\right) \qquad (5)$$

The constant "C" is the Euler constant, and its value is taken as C = 0.577215.
After performing several mathematical operations in equation (5), we obtain the following quadratic equation:

$$2aG_0 \left(\frac{1}{3} + \frac{1}{\pi^2}\right)\ell_1^2 + \left[P_1(0,t) - P_1(0,t_1)\right]\ell_1 - 2c^2 G_0 (1-C)(t-t_1) = 0 \qquad (6)$$

Thus, by accepting only the positive root of the quadratic equation (6), we obtain the following analytical formula for determining the location of the automatic valves installed on the connecting pipes at the points $x = \ell_1$ and $x = \ell_3$

$$\ell_1 = \frac{\sqrt{Z^2 + 8c^2 G_0 (1-C)(t-t_1)Z_1} - Z}{2Z_1}$$

(7)

Here, $Z = P_1(0,t) - P_1(0,t_1)$ ; $Z_1 = 2aG_0 \left(\frac{1}{3} + \frac{2}{\pi^2}\right)$

As mentioned, once the parallel gas pipelines are commissioned and the value of $\ell$, representing the step length between the connecting pipes, is known at the control center, the following analytical expression is obtained to determine the location of the automatic valves installed at the $x = \ell_3$

$$\ell_3 = \ell_1 + \ell \qquad (8)$$

As a result of the analysis conducted above, it can be noted that the technological determination of the activation time for the valves installed in the connecting pipeline should be based on the pressure values at x=0 and x= $\ell_1$ points of the damaged pipeline at t=$t_2$. Therefore, for the transmission of information to the control center and the reliable management of the gas pipeline, it is crucial to initially receive data from the sensors. To analyze and verify the data obtained at the control center, the following system equation is adopted for the creation of an algorithm based on machine learning mechanisms.

$$\begin{cases} \dfrac{P_1}{P(0,t)} \prec \varepsilon \\ \dfrac{P_1(\ell_1,0)}{P(\ell_1,t)} \succ 1 \end{cases} \quad (9)$$

In the initial stage, based on the changes in pressure values at both the starting and ending sections of the gas pipeline, the control center can record the valve activation status. Using equations (7) and (8), the location of the connecting pipelines to the right and left of the leak point is determined. Subsequently, the analysis of the data from the pressure sensors at the control center confirms condition (9). This allows for the detection of real-time operation, i.e., the time t=$t_2$. At this moment, the valves on the connecting pipeline are activated from the control center. As a result, consumers of different categories supplied by the parallel gas pipeline will be continuously provided with gas.

As seen in Table 1 and Figure 2, during the time period t=$t_1$+120 seconds, the pressure at the damaged section of x= $\ell_1$ is greater than the pressure in the undamaged pipeline (12.91×10$^4$ > 12.19×10$^4$), and at this moment, $\varepsilon$=14.58/14=1,05<1.3. Therefore, condition (9) is satisfied. This means that, for the given accident scenario and the specific gas pipeline, the control center can activate the valves installed on the connecting pipelines at $t_2$ =$t_1$+120 seconds. For the purpose of verifying the above findings and the derived empirical formula, we refer to the theoretical experimental data.

According to the data from Table 1 and the specific gas pipeline information, using equations (7) and (8), we can determine the location of the valve to be activated on the connecting pipeline, using the pressure value at the starting point of the gas pipeline at t=$t_1$+120 seconds.

$$Z = P_1(0,t) - P_1(0,t_1) = 14.58 - 13.36 = 1.22$$

$$Z_1 = 2aG_0 \left( \dfrac{1}{3} + \dfrac{2}{\pi^2} \right) = 0.1 \cdot 10 \cdot (0.333 + 0.203) = 0.536$$

$$\ell_1 = \dfrac{\sqrt{Z^2 + 8c^2 G_0 (1-C)(t-t_1) Z_1} - Z}{2 Z_1} = 0.9 \cdot 10^4 \, \text{m} ; \quad \ell_3 = 0.9 \cdot 10^4 + 1 \cdot 10^4 = 1.9 \cdot 10^4 \, \text{m}$$

To obtain theoretical experimental results, the actual distance between the connecting pipes in the specific parallel gas pipeline was accepted as $\ell_1$=1×10$^4$ m and $\ell_3$=2×10$^4$ m. Thus, the relative error between the actual value of the valve installed on the connecting pipeline of the specific parallel gas pipeline and the value determined using equation (7) is 0.1. This value of the relative error is based on the linearization of the non-linear differential equation system of non-stationary gas flow in pipeline systems [28] and the simplification of equation (1). However, this indicator is acceptable for determining the location of the connecting pipeline because the distance between the connecting pipelines installed on parallel gas pipelines is sufficiently large from both a technological and economic perspective. For a sufficiently long distance, the relative error mentioned above does not interfere with identifying the location of the connecting pipelines. For instance, if the length of the parallel gas pipeline is at least 10 km, the distance between the connecting pipes will be at least 1 km from a technological and economic standpoint [13]. Therefore, the relative error value does not hinder the identification of the location of the connecting pipes. Finally, the established formulas (7) and (8) are useful for engineering calculations and for the central control center.

After determining the location of the connecting pipes, the data obtained from the pressure sensors installed at that location are used to determine inequality (9). If the condition is met, the pneumatic valves installed at those points are closed by the control center at $t= t_2$. Thus, the consumers supplied by the parallel gas pipeline will be continuously provided with gas. To manage the aforementioned processes, we adopt the following technological system based on Figure 2. (Figure 4).

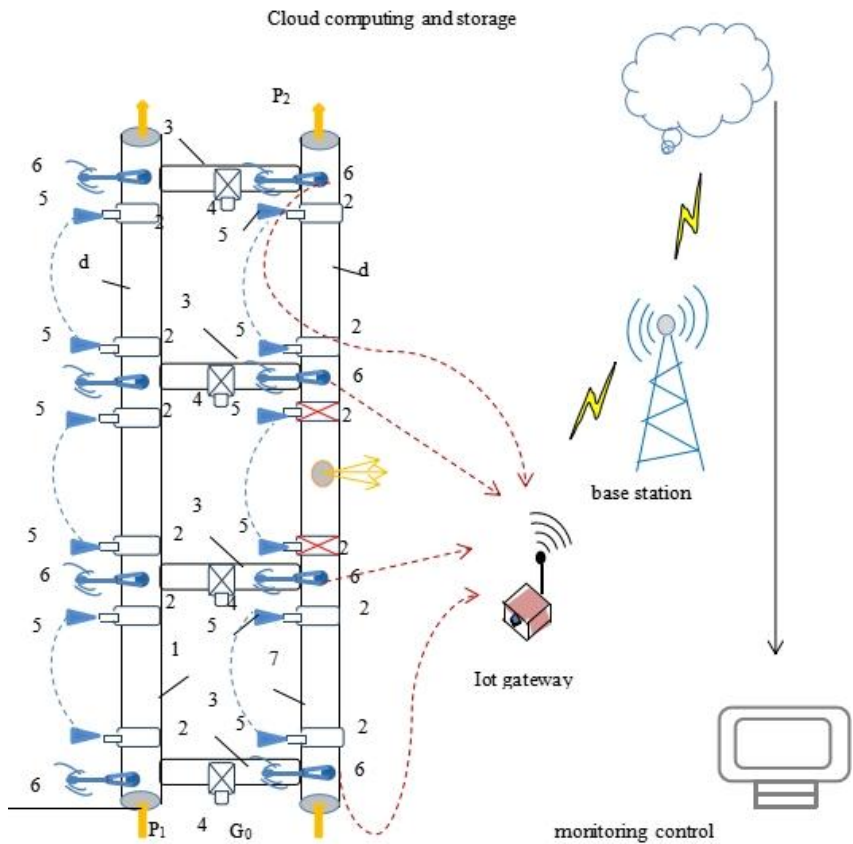

Figure 4. Technological scheme of the IoT and Cloud computing and data storage system for managing non-stationary gas flow in a parallel gas pipeline.

For monitoring the reliability of managing the non-stationary gas flow in the parallel gas pipeline, the implementation of the IoT and cloud computing and data storage system at the control center will be carried out in the following sequence, based on Figure 4.

As seen in Figure 4, when the integrity (tightness) of the 7th pipeline layout is compromised in the 8th section, the automatic valves on both sides of that section close simultaneously at time $t=t_1$ based on the control sensor's impulse. The central control center records this event based on the data from the pressure sensors located at the ends of the pipeline, $[P_1(0,t)-P_1>0$ and $P_2(L,t)-P_2<0]$, indicating that the valves are closing.

This is because after that point, the pressure at the beginning of the pipeline increases, while the pressure at the final section decreases. Such a gas dynamic process can only occur after the automatic valves have closed.

At the control center, using the equations (7) and (8), the locations of the connecting pipes on both sides of the damaged section are determined. The data from the pressure sensors installed at those locations is transferred to the Internet of Things (IoT) gateway. From the IoT gateway, the data is forwarded to the cloud storage system for storage and analysis through the base station equipment. The data received from the control center's monitoring and control devices is analyzed and calculated. The calculated parameter values are determined, and if the condition of the system equation (9) is satisfied, the valves installed on the connecting pipes (x= $\ell_1$ and x= $\ell_3$) are opened by the personnel at the control center at time t=$t_2$. Thus, consumers fed by the parallel gas pipeline will be continuously supplied with gas.

The damaged section of the parallel gas pipeline is repaired by the personnel of a special emergency brigade according to technical-operational requirements for a specific period. After the repairs, the parallel gas pipeline is put back into operation in its previous stationary mode.

## 5. CONCLUSION

As a result of the analysis of the mathematical model of non-stationary gas flow reflecting the control processes of parallel gas pipelines, an empirical algorithm based on machine learning technology has been defined for the information base of the control center. The application of the developed algorithm to the technological control system will increase the efficiency of automated decision-making processes.

The algorithms designed based on the developed schematic can be recommended for the reconstruction of existing parallel gas pipeline systems to solve engineering problems. This will enable the improvement of the central control system for the efficient operation of parallel gas pipelines.

The analysis shows that the pressure changes at the start and end sections of the parallel pipeline are directly dependent on the valve operation time. In the event of a failure, the rapid closure of the valves through the pressure sensors in the damaged sections of the parallel gas pipeline leads to the isolation of the damaged pipeline from the main network, significantly reducing gas loss to the environment and minimizing technological losses.

An empirical formula has been derived for determining the location of the remotely activated connecting valves in real-time, based on the analysis of data received from pressure sensors in the event of a failure. The use of this formula at the central control center guarantees the uninterrupted gas supply to consumers, thus enhancing the reliability and safety of the parallel gas pipeline.

A technological scheme ensuring the application of modern control systems during the reconstruction phase of existing parallel gas pipelines has been proposed to guarantee the reliability and safety of parallel gas pipelines. In order to advance this technology and apply it to operational gas pipelines, ongoing scientific and practical research in this area is essential. Consequently, it is recommended to implement an efficient management system utilizing IoT and cloud technologies in parallel gas networks.

## REFERENCES.